\documentclass[12pt]{amsart}
\usepackage[top=1in,left=1in,bottom=1in,right=1in]{geometry}
\usepackage{amsmath}
\usepackage{amssymb}
\usepackage{amsthm}
\usepackage{tikz}
\usetikzlibrary{arrows}
\usepackage{wrapfig}
\usepackage{graphicx}
\usepackage{epstopdf}
\usepackage{tikz-cd}
\usepackage{url}
\usepackage{scrextend}
\usepackage{colordvi}
\usepackage{color}
\usepackage{verbatim}

\usepackage{amsmath}
\usepackage{amsfonts}
\usepackage{colordvi}
\usepackage{color}
\usepackage{amsthm}
\usepackage{url}
\usepackage{graphicx}
\usepackage{epstopdf}
\usepackage{amssymb}

\newtheorem{theorem}{Theorem}[section]
\newtheorem{lemma}[theorem]{Lemma}
\newtheorem{conjecture}[theorem]{Conjecture}

\newtheorem{question}[theorem]{Question}

\newtheorem{proposition}[theorem]{Proposition}

\theoremstyle{definition}
\newtheorem{definition}[theorem]{Definition}
\theoremstyle{remark}
\newtheorem{remark}[theorem]{Remark}
\newtheorem{example}[theorem]{Example}

\newcommand{\A}{\mathcal A}

\newcommand{\C}{\mathcal C}
\newcommand{\D}{\mathcal D}

\renewcommand{\S}{\mathcal S}

\newcommand{\U}{\mathcal U}

\newcommand{\X}{\mathcal X}

\newcommand{\R}{\mathbb R}

\DeclareMathOperator{\conv}{conv}

\DeclareMathOperator{\cl}{cl}

\newcommand{\mxl}[1]{\mathbf{#1}}

\DeclareMathOperator{\nddim}{nddim}
\DeclareMathOperator{\odim}{odim}
\DeclareMathOperator{\cdim}{cdim}
\DeclareMathOperator{\code}{code}
\DeclareMathOperator{\interior}{int}

\newcommand{\od}{:=}

\renewcommand{\S}{\mathcal S}
\renewcommand{\mxl}[1]{#1}

\begin{document}
\title{Open, Closed, and Non-Degenerate Embedding Dimensions of Neural Codes}

\author{R. Amzi Jeffs}

\begin{abstract}
We study the open, closed, and non-degenerate embedding dimensions of neural codes, which are the smallest respective dimensions in which one can find a realization of a code consisting of convex sets that are open, closed, or non-degenerate in a sense defined by Cruz, Giusti, Itskov, and Kronholm. For a given code $\C$ we define the embedding dimension vector to be the triple $(a,b,c)$ consisting of these embedding dimensions. Existing results guarantee that $\max\{a,b\} \le c$, and we show that when any of these dimensions is at least 2 this is the only restriction on such vectors. Specifically, for every triple $(a,b,c)$ with $2\le \min \{a,b\}$ and $\max\{a,b\}\le c\le \infty $ we construct a code $\C_{(a,b,c)}$ whose embedding dimension vector is exactly $(a,b,c)$ (where an embedding dimension is $\infty$ if there is no realization of the corresponding type). 

Our constructions combine two existing tools in the convex neural codes literature: sunflowers of convex open sets, and rigid structures, the latter of which was recently defined in work of Chan, Johnston, Lent, Ruys de Perez, and Shiu. Our constructions provide the first examples of codes whose closed embedding dimension is larger than their open embedding dimension, but still finite. 
\end{abstract}

\thanks{Jeffs' work on this paper was supported by the National Science Foundation through grants DGE-1761124 and Award No. 2103206.}
\date{\today.Department of Mathematics, Carnegie Mellon University. Wean Hall, 5000 Forbes Ave, Pittsburgh, PA 15213.\\
This version of the article has been accepted for publication, after peer review but is not the Version of Record and does not reflect post-acceptance improvements, or any corrections. The Version of Record is available online at: \url{https://doi.org/10.1007/s00454-023-00512-1}.}
\maketitle

\section{Introduction}\label{sec:intro}

In 2013 Curto, Itskov, Veliz-Cuba, and Youngs \cite{neuralring13} initiated the study of convex neural codes, which are the combinatorial codes $\C\subseteq 2^{[n]}$ that record the intersection and covering relations among $n$ convex open sets in Euclidean space (see Definition \ref{def:codeofrealization} below). Their motivation arose from neuroscience, namely the study of place cells, which are hippocampal neurons that fire when an animal is in a particular region of its environment. Place cells can be thought of as encoding a cognitive map of the animal's environment (see \cite{okeefe}), and the study of convex neural codes seeks to understand how well this cognitive map can capture the geometry and topology of the environment. 

Mathematical research on convex codes has blossomed since 2013. An efficient characterization of convex codes is unfortunately out of the question---recent work in \cite{matroids} shows that recognizing convex codes is $\exists \R$-hard. Nevertheless, researchers have used techniques from algebra \cite{signatures, grobner, polarization}, discrete geometry \cite{rigidstructures, openclosed, sunflowers, CUR, obstructions}, and topology \cite{undecidability, local15} to analyze many interesting families of codes and develop frameworks in which to test whether or not a code is convex. Some works also study codes that arise from ``good covers" \cite{undecidability, local15}, collections of closed convex sets \cite{openclosed, 5neurons}, and ``non-degenerate" collections of convex sets \cite{rigidstructures, openclosed}. 

Let us begin by recalling some fundamental definitions. Our combinatorial objects of study are codes, which are just subsets of the power set $2^{[n]}$, where $[n]\od \{1,2,\ldots, n\}$. 

\begin{definition}\label{def:code}
A collection $\C\subseteq 2^{[n]}$ is called a \emph{code}. The elements of $\C$ are called \emph{codewords}. If $\C\subseteq 2^\sigma$ for some $\sigma\subseteq [n]$, we say that $\sigma$ is a \emph{base set} for $\C$. \end{definition}

We will adopt the convention that $\emptyset$ is a codeword in every code, which is typical in the study of convex neural codes. Codes can be used to record information about how sets in a collection intersect and cover one another, as follows. 

\begin{definition}\label{def:codeofrealization}
Let $\U = \{U_1, U_2, \ldots, U_n\}$ be a collection of sets in $\R^d$. The \emph{code} of $\U$ is \begin{align*}
\code(\U) & \od \{\sigma\subseteq [n] \mid \text{There exists $p\in \R^d$ such that $p\in U_i$ if and only if $i\in \sigma$}\}\\
& \, \, = \bigg\{ \sigma\subseteq [n] \,\bigg|\, \bigcap_{i\in\sigma} U_i \setminus \bigcup_{j\in[n]\setminus\sigma} U_j \neq \emptyset \bigg\} 
\end{align*}
where the empty-indexed intersection is equal to $\R^d$ by convention. We say that $\U$ is a \emph{realization} of $\code(\U)$. 
\end{definition}

In words, we obtain $\code(\U)$ by labeling every point $p\in \R^d$ by the indices $i\in[n]$ for which $p\in U_i$, and then collecting all the labels obtained in this way. If $\U$ consists of convex open sets we say that $\U$ is a \emph{convex open} realization of $\code(\U)$. We may similarly define \emph{convex closed} realizations, and we typically use the notation $\X = \{X_1, \ldots, X_n\}$ for collections of closed convex sets. Every realization in this paper will consist of convex sets, so we will usually drop the adjective ``convex." Note that our convention that $\emptyset$ is always a codeword amounts to the requirement that a realization $\U$ does not cover $\R^d$---in particular, by intersecting with a sufficiently large closed or open ball we may assume all of our realizations are bounded. 

\begin{example}\label{ex:realization}
Consider the code \[
\C = \{123, 12, 13, 24, 34, 1, 2, 3, 4, \emptyset\}.
\]
Figure \ref{fig:realization} shows an open realization $\U = \{U_1, U_2, U_3, U_4\}$ of $\C$ in $\R^2$. One could also regard this figure as an illustration of a closed realization: replacing each set by its closure does not change the realized code.  
\begin{figure}[h]
\[\includegraphics{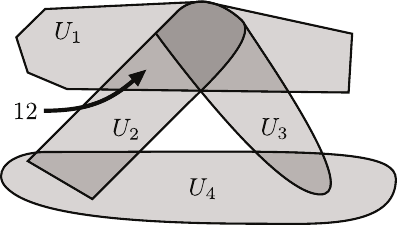}\]
\caption{An (open or closed) realization of $\C$ in $\R^2$, with an arrow pointing to the region where the codeword $12$ arises.}\label{fig:realization}
\end{figure}
\end{example}

\begin{remark}
We will always illustrate open convex sets with a solid border. To avoid confusing these illustrations with closed convex sets, our captions will always specify whether we are regarding the illustrated sets as closed or open. 
\end{remark}

In Example \ref{ex:realization}, it was convenient that we could regard our realization as either closed or open without changing the realized code. This is not always the case—for example, in an open realization we may have disjoint sets which share boundary points, so that replacing them by their closures changes the code that they realize. In fact, this may be the case in \emph{every} open realization of a code: \cite{openclosed} gives an example of a code in which every open realization is forced to include disjoint sets that share boundary points. Motivated by this difficulty, \cite{openclosed} introduced a notion of non-degeneracy for realizations, which places technical geometric and topological criteria on a realization in such a way that replacing sets by their interiors or closures preserves  the realized code. Recently, \cite{rigidstructures} proved that non-degenerate realizations are \emph{exactly} those for which replacing sets by interiors or closures does not affect the realized code---we will take this as the definition of non-degeneracy.  

\begin{definition}\label{def:nondeg}
A collection $\U = \{U_1,\ldots, U_n\}$ of convex open sets is called \emph{non-degenerate} if the collection $\X = \{X_1, \ldots, X_n\}$ with $X_i\od \cl(U_i)$ has the property that $\code(\X) = \code(\U)$. Symmetrically, a collection of closed convex sets $\X = \{X_1, \ldots, X_n\}$ is called \emph{non-degenerate} if the collection $\U = \{U_1, \ldots, U_n\}$ with $U_i \od \interior(X_i)$ has the property that $\code(\U) = \code(\X)$.  
\end{definition}

It is of particular interest to determine the smallest dimension in which a code has an open, closed, or non-degenerate realization. These minimum dimensions are referred to as \emph{embedding dimensions} of a code.

\begin{definition}\label{def:odim}\label{def:cdim}\label{def:nddim}
Let $\C$ be a code. The \emph{open}, \emph{closed}, and \emph{non-degenerate} embedding dimensions of $\C$ are the following quantities, respectively:\begin{align*}
\odim(\C) \od &\min\{d\mid \text{$\C$ has an open convex realization in $\R^d$}\},\\
\cdim(\C) \od &\min\{d\mid \text{$\C$ has a closed convex realization in $\R^d$}\}, \text{and}\\
\nddim(\C) \od &\min\{d\mid \text{$\C$ has a non-degenerate (open or closed) convex realization in $\R^d$}\}.
\end{align*}
Above, the minimum over the empty set is equal to $\infty$ by convention. The \emph{embedding dimension vector} of $\C$ is the 3-tuple \[
\big(\odim(\C), \cdim(\C), \nddim(\C)\big).
\]
\end{definition}

For a fixed code $\C\subseteq 2^{[n]}$, what can we say about the embedding dimensions of $\C$ and their relationships to one another? Determining these dimensions exactly is often an infeasible task, but it is possible in some specific cases, and sometimes one may obtain general bounds that are interesting even if they are not exact. For example, when $\C$ is \emph{intersection complete} (i.e. the intersection of any two codewords is again a codeword),  \cite{openclosed} showed that $\nddim(\C) \le \max\{2, m\}$ where $m+1$ is the number of inclusion-maximal codewords in $\C$, and \cite{embeddingphenomena} showed that $\cdim(\C)\le 2d+1$ if every codeword has size $d+1$ or less. 

As a very basic start, every non-degenerate realization can be regarded as a closed or open realization, so we have the following:

\begin{proposition}
If $\C$ has embedding dimension vector $(a,b,c)$, then $\max\{a,b\} \le c$.
\end{proposition}

It is natural to ask whether we can guarantee stricter relationships between the various embedding dimensions of a code. As we will see in Theorem \ref{thm:main}, the answer in general is ``no." However, in some special cases, the answer is yes. For example, if $\C$ is a simplicial complex then $\cdim(\C) = \odim(\C) = \nddim(\C)$ (see \cite[Theorem 1.4]{embeddingphenomena}). If $\C$ is intersection complete, then $\cdim(\C) \le \odim(\C) = \nddim(\C)$ (see \cite[Lemma 2.2.4 and Theorem 2.2.7]{amziphdthesis}).

One final special case is when any of the embedding dimensions is equal to 1. In this case, all embedding dimensions must be equal to 1. This fact was first posed as a conjecture in \cite[Conjecture 3.4]{rigidstructures}, and we prove it below using ideas based on discussions with the authors.

\begin{theorem}\label{thm:1dim}
Let $\C$ be a code. Then the following are equivalent:\begin{itemize}
\item[(i)] $\odim(\C) = 1$,
\item[(ii)] $\cdim(\C) = 1$, and
\item[(iii)] $\nddim(\C) = 1$.
\end{itemize}
\end{theorem}
\begin{proof}
It will suffice to show that any open or closed realization of $\C$ by intervals in $\R^1$ can be made non-degenerate. In any realization by (open or closed) intervals, we may assume without loss of generality that every point in $\R^1$ is either a left endpoint of some intervals in our realization, or a right endpoint of some intervals in our realization, but not both simultaneously. To guarantee this, simply insert a closed unit interval $[a,b]$ at any point $p\in \R^1$ that is simultaneously a left and right endpoint. If our realization is open, we modify it so that $a$ is a right endpoint of all intervals that $p$ was a right endpoint of, and $b$ is a left endpoint of all intervals that $p$ was a left endpoint of. If our realization is closed, we do the opposite: intervals whose left endpoints were equal to $p$ now have left endpoint $a$, while those with right endpoint $p$ now have right endpoint $b$. 

We claim that such a realization is necessarily non-degenerate. Observe that if some codeword $c$ arises at a point $p$, then the same codeword arises at every point in a small closed interval with one of its endpoints equal to $p$. Thus replacing our intervals by their interiors or closures does not change the realized code, and the realization is non-degenerate.
\end{proof}

Beyond the 1-dimensional case, the only relationship that we can guarantee between embedding dimensions is that the open and closed embedding dimensions are no larger than the non-degenerate embedding dimension. The following theorem captures this fact formally.

\begin{theorem}\label{thm:main}
Let $2\le a,b,c \le \infty$ and suppose that $\max\{a,b\} \le c$.
Then there exists a code $\C_{(a,b,c)}$ with embedding dimension vector $(a,b,c)$.
\end{theorem}

Rather than construct all $\C_{(a,b,c)}$ directly, we will reduce to three cases, from which one can build any $\C_{(a,b,c)}$. For this reduction, we observe that any two codes $\C$ and $\D$ may be relabeled so that they have disjoint base sets, and then combined to yield a code whose embedding dimension vector is the component-wise maximum of the original embedding dimension vectors. 

\begin{proposition}\label{prop:max}
Let $\C$ and $\D$ be codes on disjoint base sets, with respective embedding dimension vectors $(a_1, b_1, c_1)$ and $(a_2, b_2, c_2)$. Then the embedding dimension vector of $\C\cup \D$ is \[
\big(\max\{a_1, a_2\}, \max\{b_1, b_2\}, \max\{c_1, c_2\}\big).
\]
\end{proposition}
\begin{proof}
Any realization of $\C\cup\D$ yields a realization of $\C$ by deleting the sets indexed by the base set of $\D$, and vice versa. Deleting sets preserves openness, closedness, or non-degeneracy of a realization, so the various embedding dimensions of $\C\cup\D$ provide an upper bound on the corresponding embedding dimensions of $\C$ and $\D$. Conversely, any pair of (open, closed, or non-degenerate) realizations of $\C$ and $\D$ in the same dimension yields a corresponding realization of $\C\cup \D$ by placing the two realizations sufficiently far apart. This proves the result.
\end{proof}

By Proposition \ref{prop:max}, to prove Theorem \ref{thm:main} it will suffice to exhibit codes with embedding dimension vectors $(d,2,d)$, $(2,d,d)$, and $(2,2,d)$ for all choices of $2\le d\le \infty$. We treat the respective cases for finite $d$ in Sections \ref{sec:(d,2,d)}, \ref{sec:(2,d,d)}, and \ref{sec:(2,2,d)}. The cases with $d=\infty$ are all treated in Section \ref{sec:infty}.

Our most technical result is the construction and analysis of the code $\C_{(2,d,d)}$---in particular, proving that $\C_{(2,d,d)}$ has a non-degenerate realization in $\R^d$ requires several pages of careful work (see Proposition \ref{prop:(2,d,d)nondeg}). 

Our constructions primarily make use of two existing tools. First, in \cite{embeddingphenomena, sunflowers} we studied ``sunflowers" of convex open sets, obtaining examples of codes with large open embedding dimension and small closed embedding dimension. Second, Chan, Johnston, Lent, Ruys de Perez, and Shiu introduced ``rigid structures" in \cite{rigidstructures}. Their results guarantee that sets in a closed realization must have a union which is convex under certain conditions---and importantly, their results do not hold for open realizations. 

Informally, sunflowers guarantee structure in open (but not closed) realizations, while rigid structures provide the opposite. By combining these tools in various ways we are able to obtain all the desired codes $\C_{(a,b,c)}$. In the interest of concision, we do not explain sunflowers or rigid structures in full generality. Instead, we state versions of these results that suffice in our context, and provide citations for a more general presentation. 



\section{Constructing the Codes $\C_{(d,2,d)}$ with $d< \infty$}\label{sec:(d,2,d)}

We begin by constructing the codes $\C_{(d,2,d)}$ for all finite $d$. In fact, there is an existing family in the convex codes literature that suffices: sunflower codes. We first introduced and studied these codes in \cite[Definition 5.3]{embeddingphenomena}, where we were primarily concerned with their open embedding dimensions. Below we review the definition of these codes, illustrate a few small examples, and provide citations for the results implying that they have the appropriate embedding dimension vector.

\begin{definition}[See also \cite{embeddingphenomena} and  \cite{amziphdthesis}]\label{def:Sd}
For $2\le d < \infty$, define $\S_d\subseteq 2^{[d+1]}$ to be the code consisting of the following codewords: $[d]$, all singleton sets, all pairs $\{i, d+1\}$ for $i\in [d]$, and the empty set. 
\end{definition}

The code $\S_d$ has two salient geometric features. First, in any realization the first $d$ sets must form a ``sunflower" in the sense that their various pairwise intersections must be the same, and must be nonempty. Second, the $(d+1)$-st set intersects all the other sets, but not their common intersection. It turns out that such an arrangement is only possible to achieve with convex open sets in dimension at least $d$. For a full discussion of this fact, see \cite[Section 5]{embeddingphenomena}. Below we illustrate a few of these codes to provide intuition.

\begin{example}
The sunflower codes $\S_d$ for $d = 2, 3, 4$ are listed below:\begin{align*}
\S_2 = & \{\mxl{12}, \mxl{13}, \mxl{23}, 1,2,3, \emptyset\}, \\
\S_3 = & \{\mxl{123}, \mxl{14}, \mxl{24}, \mxl{34}, 1,2,3,4, \emptyset\}, \text{and} \\
\S_4 = & \{\mxl{1234}, \mxl{15}, \mxl{25}, \mxl{35}, \mxl{45}, 1,2,3,4,5 , \emptyset\}.
\end{align*}
Figure \ref{fig:sunflowers} shows open realizations for $\S_2$ and $\S_3$ in $\R^2$ and $\R^3$ respectively. The code $\S_4$ does not have an open realization in $\R^3$, so we illustrate a closed realization in $\R^2$. 
\begin{figure}[h]
\[
\includegraphics{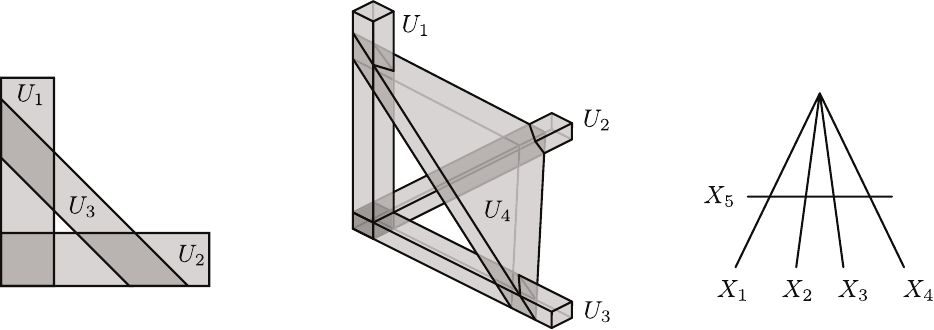}
\]
\caption{Open realizations of $\S_2$ and $\S_3$, and a closed realization of $\S_4$.}\label{fig:sunflowers}
\end{figure}
\end{example}

We conclude this section by formally observing that $\S_d$ has the desired embedding dimension vector $(d,2,d)$. We provide several citations rather than a detailed proof since the results characterizing the embedding dimensions of $\S_d$ are already established.

\begin{proposition}\label{prop:(d,2,d)}
Let $2\le d < \infty$. Then the code $\S_d\subseteq 2^{[d+1]}$ has embedding dimension vector equal to $(d,2,d)$. 
\end{proposition}
\begin{proof}
The results \cite[Theorem 5.2.2 and Proposition 5.2.3]{amziphdthesis} tell us that $\odim(\S_d) = d$ and $\cdim(\S_d) = 2$. The code $\S_d$ is intersection complete, and \cite[Lemma 2.2.4]{amziphdthesis} guarantees that open and non-degenerate embedding dimension are equal for intersection complete codes. Thus $\nddim(\S_d) = \odim(\S_d) = d$, proving the result.
\end{proof}

\section{Constructing the Codes $\C_{(2,d,d)}$ with $d< \infty$}\label{sec:(2,d,d)}

The code $\C_{(2,d,d)}$ will have a base set of size $4d$. Rather than simply use the integers $\{1,2,\ldots, 4d\}$, we use four types of labeled symbol, with $d$-many of each. Below, we let \[\alpha_{[d]} \od \{\alpha_1, \alpha_2, \ldots, \alpha_d\}\] where each $\alpha_i$ is a formal symbol in our base set. For any $\sigma\subseteq [d]$, we let $\alpha_{\sigma} = \{\alpha_i\mid i\in \sigma\}$, and in specific examples we will sometimes omit braces on $\sigma$---so for example, $\alpha_{23} = \{\alpha_2, \alpha_3\}$. The sets $\beta_{[d]}$, $\gamma_{[d]}$ and $\delta_{[d]}$ are defined analogously.

This notation has two distinct advantages. First, it streamlines the indexing in the results below, so that we need not deal with cumbersome offset factors in our base set indices. Second, it highlights that each type of base set element plays a different role in the code. Before commenting on these various roles, we provide a formal definition. 

\begin{definition}\label{def:(2,d,d)}
Let $2\le d < \infty$, and define $\C_{(2,d,d)}$ be the code on the base set $\alpha_{[d]}\cup\beta_{[d]}\cup\gamma_{[d]}\cup\delta_{[d]}$ which has the following nonempty codewords:\begin{itemize}
\item[(i)] $\{\beta_i, \gamma_i\}$ for $i\in [d]$,
\item[(ii)] $\{\beta_i,\beta_{i+1}, \gamma_i\}$ for $i\in[d-1]$,
\item[(iii)] $\{\beta_i,\beta_{i+1}\}$ for $i\in[d-1]$,
\item[(iv)] $\{\beta_i,\beta_{i+1}, \gamma_{i+1}\}$ for $i\in[d-1]$,
\item[(v)] $\{\alpha_i, \beta_i, \gamma_i,\delta_i \}$ for $i\in [d]$,
\item[(vi)] $\{\alpha_i,\delta_i\}$ for $i\in [d]$,
\item[(vii)] $\alpha_{[d]}\cup \delta_{[d]}$,
\item[(viii)] $\delta_{[d]}$.
\end{itemize}
\end{definition}

Informally, the base set elements and codewords above each play the following roles. The various $\alpha_i$ are defined so that the various $U_{\alpha_i}$ in any realization of $\C_{(2,d,d)}$ form  a sunflower---recall the commentary following Definition \ref{def:Sd}---thanks to the codewords of type (vi) and (vii). The codewords of types (i)-(iv) guarantee that the various $\beta_i$ and $\delta_i$ form a ``rigid structure" as defined in \cite{rigidstructures}---this means that the union of all $X_{\beta_i}$ and $X_{\gamma_i}$ in a closed realization must be convex (see Lemma \ref{lem:rigid} below). Moreover, the codewords of type (v) force this rigid structure to intersect the various sunflower petals $U_{\alpha_i}$. Finally, the various $\delta_i$ have essentially the same behavior as the $\alpha_i$, with the exception of the codeword (viii), which ties the structure of $\C_{(2,d,d)}$ to the structure of the code $\A_d$ from \cite[Theorem 3.7]{nonmonotonicity}, and is key to forcing the closed embedding dimension of $\C_{(2,d,d)}$ to be large. 

Let us start our analysis of $\C_{(2,d,d)}$ concretely,  by forming an open realization in $\R^2$. 

\begin{proposition}\label{prop:(2,d,d)open}
The code $\C_{(2,d,d)}$ has an open realization in $\R^2$. 
\end{proposition}
\begin{proof}
We begin by describing the sets $U_{\alpha_i}$ and $U_{\delta_i}$ for $i\in [d]$. Let $P$ be a regular $(d+1)$-gon in $\R^2$ with center at the origin, inscribed in a second regular $(d+1)$-gon $P'$, which is rotated by an angle of $\pi/(d+1)$ so that the vertices of $P$ meet the midpoints of the edges of $P'$. Observe that $\interior(P'\setminus P)$ consists of $d+1$ disjoint open ``flaps" arranged sequentially along the edges of $P$. Label these flaps as $F_1, F_2, \ldots, F_{d+1}$. For each $i\in[d]$ define \begin{align*}
U_{\alpha_i} &= \interior(P\cup F_i), \text{and}\\
U_{\delta_i} &= \interior(P\cup F_i \cup F_{d+1}).
\end{align*}
In words, $U_{\alpha_i}$ is $P$ plus the $i$-th flap, while $U_{\delta_i}$ is $P$ plus the $i$-th flap and the $(d+1)$-st flap. Observe that the nonempty codewords arising in this arrangement are exactly those of types (vi), (vii), and (viii) in Definition \ref{def:(2,d,d)}. Indeed, the codeword $\{\alpha_i, \delta_i\}$ arises in $F_i$, the codeword $\alpha_{[d]}\cup \delta_{[d]}$ arises in the interior of $P$, and the codeword $\delta_{[d]}$ arises inside $F_{d+1}$. 

We can now define the sets $U_{\beta_i}$ and $U_{\gamma_i}$ in our realization. For each $i\in[d+1]$, let $L_i$ be a line that is parallel to the $i$-th edge of $P$, moved a small distance away from $P$ but still intersecting the flap $F_i$. Let $L_i'$ be a second copy of $L_i$ moved twice as far from $P$ as $L_i$, and note that $L_i'$ does not necessarily pass through $F_i$. For $2 \le i \le d$, label points in the intersections of the various $L_i$ and $L_i'$ as follows:\begin{align*}
p_i &= L_{i-1}\cap L_{i} & \quad q_i &= L_{i-1}'\cap L_{i}'\\
r_i &= L_{i-1}\cap L_{i}' & \quad s_i &= L_{i-1}'\cap L_{i}\\
t_i & = (q_i+r_i)/2 &\quad u_i &= (q_i+s_i)/2. 
\end{align*} These points are shown in Figure \ref{fig:pqrstu}. Moreover, we define the following edge cases:\begin{align*}
p_1 & = s_1 = L_1\cap L_{d+1}' &\quad q_1 &=t_1 = L_{1}' \cap L_{d+1}'\\
r_{d+1} &=p_{d+1} = L_{d}\cap L_{d+1}' &\quad q_{d+1} &= u_{d+1} = L_d'\cap L_{d+1}'. 
\end{align*} 
\begin{figure}[h]
\[
\includegraphics{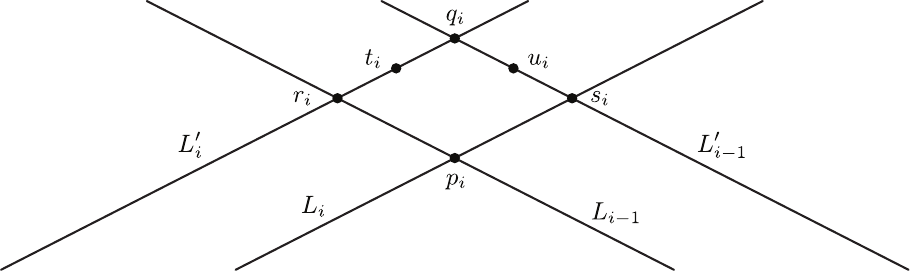}\]
\caption{Some of the points used to construct $U_{\beta_i}$ and $U_{\gamma_i}$. The polygons $P$ and $P'$ are not pictured.}\label{fig:pqrstu}
\end{figure}
Now, with these points labeled, for $i\in[d]$ we define:\begin{align*}
U_{\beta_i} &= \interior(\conv\{s_{i}, r_{i+1}, q_{i+1}, q_{i}\}), \text{and} \\
U_{\gamma_i} &= \interior(\conv\{p_{i}, p_{i+1}, u_{i+1}, t_{i}\}).
\end{align*}
We claim that this completes our open realization of $\C_{(2,d,d)}$. The $U_{\beta_i}$ and $U_{\gamma_i}$ do not fully cover any of the regions giving rise to the codewords of types (vi)-(viii) that we described previously, so it suffices to show that the codewords arising inside the various $U_{\beta_i}$ and $U_{\gamma_i}$ are exactly those of types (i)-(v) in Definition \ref{def:(2,d,d)}.

Note that the union of the various $U_{\beta_i}$ and $U_{\gamma_i}$ is a bent, thickened line segment which wraps around the first $d$ edges of $P$. For $i\in[d-1]$, the codewords arising near the $i$-th joint in this bent region are exactly $\{\beta_i, \gamma_i\}$, $\{\beta_i,\beta_{i+1}, \gamma_i\}$, $\{\beta_i,\beta_{i+1}\}$, $\{\beta_i, \beta_{i+1}, \gamma_{i+1}\}$, $\{\beta_{i+1}, \gamma_{i+1}\}$ in sequence, as illustrated in Figure \ref{fig:joint}. 

\begin{figure}[h]
\[\includegraphics{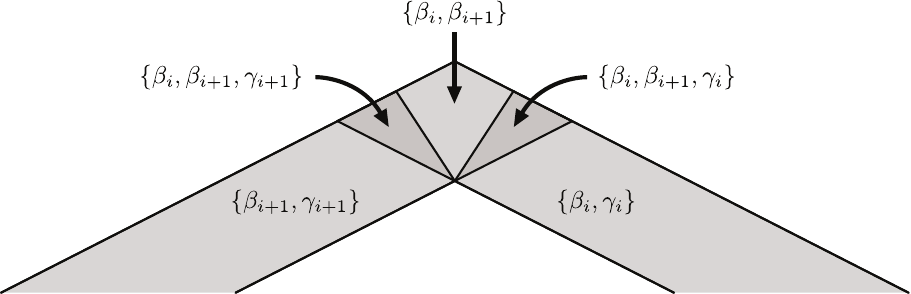}\]
\caption{The codewords arising at the $i$-th ``joint" where $L_i$ and $L_{i+1}$ meet. Note that this open realization is degenerate: the disjoint sets $U_{\gamma_i}$ and $U_{\gamma_{i+1}}$ share the boundary point $p_{i+1} = L_i\cap L_{i+1}$. }\label{fig:joint}
\end{figure}

These are all codewords of type (i)-(iv), and we see that all such codewords arise for various choices of $i$. Away from the joints, the only codewords that arise are $\{\beta_i, \gamma_i\}$ and $\{\alpha_i, \beta_i, \gamma_i, \delta_i\}$, the latter arising near the midpoint of $L_i$. This accounts for the codewords of type (v), and so we have indeed realized the code $\C_{(2,d,d)}$ as desired. 
\end{proof}

\begin{example}\label{ex:(2,d,d)open}
In Figure \ref{fig:(2,d,d)open} we illustrate the open realization of $\C_{(2,d,d)}$ constructed in Proposition \ref{prop:(2,d,d)open} in the case $d=4$. Explicitly, the code we realize is \begin{align*}
\C_{(2,4,4)} = & \{\beta_1\gamma_1, \, \, \alpha_1\beta_1\gamma_1\delta_1, \, \,  \beta_{12}\gamma_1, \, \, \beta_{12}, \, \,  \beta_{12}\gamma_2, \\
& \, \,\, \beta_2\gamma_2, \, \,  \alpha_2\beta_2\gamma_2\delta_2, \, \, \beta_{23}\gamma_2, \, \,  \beta_{23}, \, \,  \beta_{23}\gamma_3, \\
& \, \, \,
\beta_3\gamma_3, \, \,  \alpha_3\beta_3\gamma_3\delta_3, \, \,  \beta_{34}\gamma_3, \, \,  \beta_{34}, \, \,  \beta_{34}\gamma_4,\\
&\,\,\, \beta_4\gamma_4, \,\, \alpha_4\beta_4\gamma_4\delta_4, \\
& \, \,\, \alpha_1\delta_1, \, \,  \alpha_2\delta_2, \, \,  \alpha_3\delta_3, \, \,  \alpha_4\delta_4, \, \,  \delta_{1234}, \, \,  \alpha_{1234}\delta_{1234}, \, \, \emptyset\}.
\end{align*}
In the first four lines above we have written the codewords that appear in the bent region around the outside of the central pentagon in the order that they appear. In the final line we have written the codewords that appear in the five ``flaps" around the pentagon and the codeword $\alpha_{1234}\delta_{1234}$ that appears in the pentagon itself.  
\begin{figure}
\[\includegraphics{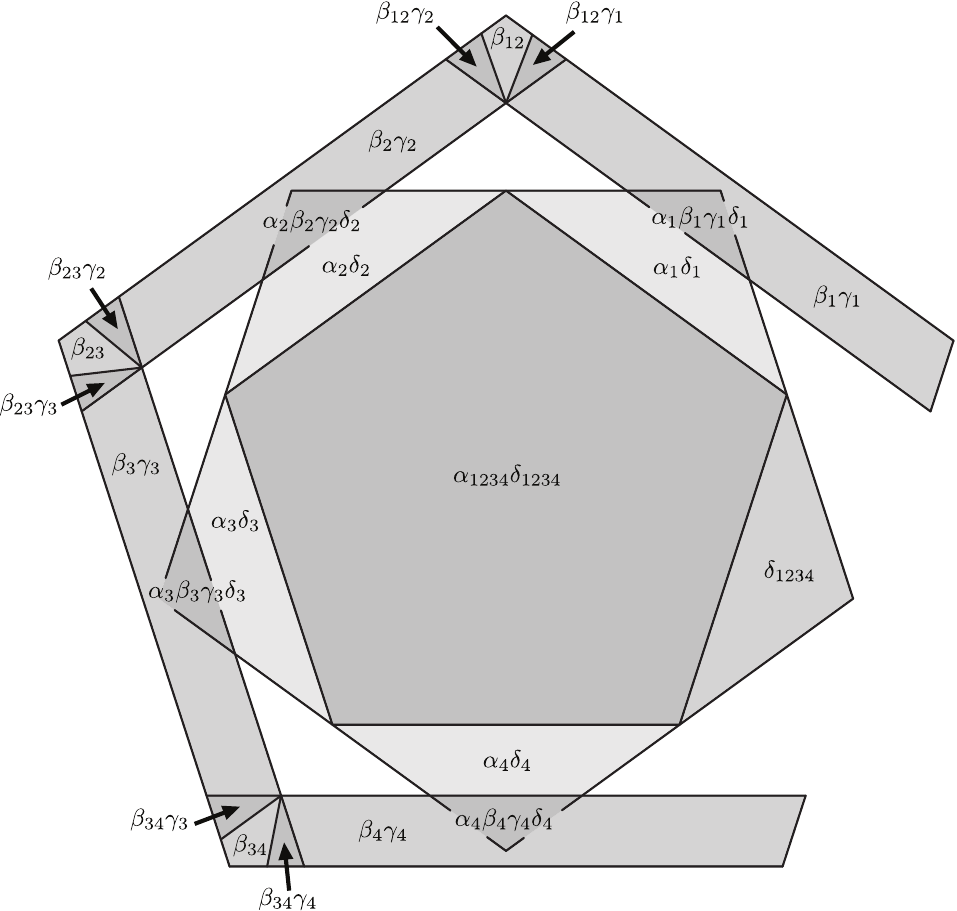}\]
\caption{An open realization of $\C_{(2,4,4)}$ in $\R^2$ as constructed in the proof of Proposition \ref{prop:(2,d,d)open}.}\label{fig:(2,d,d)open}
\end{figure}
\end{example}

With this construction achieved, we can proceed to construct a non-degenerate open (and hence also closed) realization of $\C_{(2,d,d)}$ in $\R^d$. This construction is the most technical result in the paper, and requires us to carefully manipulate a variety of inequalities that define the sets in our realization. However, the broad intuition for this construction is not too complex: we thicken the coordinate axes into open cubical prisms to form the various $U_{\alpha_i}$ and $U_{\delta_i}$, and we form the various $U_{\beta_i}$ and $U_{\gamma_i}$ by sequentially slicing through a thickened simplex in the positive orthant which lies far form the origin. This construction is illustrated for the case $d=3$ in Example \ref{ex:(2,3,3)} below.

\begin{proposition}\label{prop:(2,d,d)nondeg}
The code $\C_{(2,d,d)}$ has a non-degenerate realization in $\R^d$. 
\end{proposition}
\begin{proof}
For $i\in[d]$, we start by defining \begin{align*}
U_{\alpha_i} &= \bigg \{\mathbf x\in \R^d \, \bigg|\,  0 < x_j < 1 \text{ for $j \neq i$ and } \sum_{j\in [d]} x_j > 1\bigg \}, \text{and}\\
U_{\delta_i} & = \left \{\mathbf x\in \R^d \mid 0 < x_j < 1 \text{ for $j \neq i$ and } x_i > 0\right \}. 
\end{align*}
In words, $U_{\delta_i}$ is obtained from an open unit hypercube in the positive orthant by extending it infinitely in the $i$-th coordinate direction. The set $U_{\alpha_i}$ is the subset of $U_{\delta_i}$ in which the sum of all coordinates is larger than one---that is, $U_{\alpha_i}$ is obtained from $\U_{\delta_i}$ by slicing off the simplex in which the sum of coordinates is one or less, which lies in the corner of the positive orthant. 

We claim that the nonempty codewords that arise among these sets are exactly those of types (vi), (vii), and (viii) in Definition \ref{def:(2,d,d)}. Note that outside of the open unit hypercube in which each coordinate is between zero and one, the only nonempty codewords appearing are $\{\alpha_i, \delta_i\}$ for $i\in [d]$. This follows from the fact that $U_{\alpha_i}$ and $U_{\gamma_i}$ are the same outside the hypercube, and $U_{\alpha_i}$ does not meet $U_{\alpha_j}$ outside the hypercube when $j\neq i$. Inside the hypercube there are two regions. Where the sum of coordinates is larger than one, all $U_{\alpha_i}$ and $U_{\delta_i}$ appear, giving rise to the codeword $\alpha_{[d]}\cup \delta_{[d]}$. Where the sum of coordinates is one or less, only the various $U_{\delta_i}$ appear, giving rise to the codeword $\delta_{[d]}$. Thus all codewords of types (vi), (vii), and (viii) appear, and no others.

Now let us define the $U_{\beta_i}$ and $U_{\gamma_i}$. Let $\ell$ denote the linear functional given by $\ell(\mathbf x)  = \sum_{i\in[d]} i x_i$. Let $C$ be the open convex region in the positive orthant consisting of the points $\mathbf x$ so that the sum of the coordinates of $\mathbf x$ is between $5d^2$ and $5d^2+1$.
Now, for $i\in[d]$ define \begin{align*}
U_{\beta_i} &= \{\mathbf x \in C\mid 5id^2-4d^2 < \ell(\mathbf x) < 5id^2 + 4d^2\}, \text{and}\\
U_{\gamma_i} &= \{\mathbf x \in C\mid 5id^2-2d^2 < \ell(\mathbf x) < 5id^2 + 2d^2\}.
\end{align*}
We aim to show the addition of these sets to our realization gives rise to exactly the codewords of type (i)-(v) from Definition \ref{def:(2,d,d)}. First, let us determine the codewords that arise from these sets independent of the $U_{\alpha_i}$ and $U_{\delta_i}$. If $\mathbf x\in C$, then observe that the value of $\ell(\mathbf x)$ completely determines which codeword arises at $\mathbf x$:
\begin{itemize}
\item $\{\beta_i,\gamma_i\}$ for $i\in [d]$ arises when $5id^2 - d^2 \le \ell(\mathbf x) \le 5id^2 + d^2$, 
\item $\{\beta_i,\beta_{i+1}, \gamma_i\}$ for $i\in[d-1]$ arises when $5id^2 + d^2 < \ell(\mathbf x) < 5id^2 + 2d^2$,
\item $\{\beta_i,\beta_{i+1}\}$ for $i\in[d-1]$ arises when $5id^2 + 2d^2 \le\ell(\mathbf x) \le 5id^2 + 3d^2$, and
\item $\{\beta_i,\beta_{i+1},\gamma_{i+1}\}$ for $i\in[d-1]$ arises when $5id^2 + 3d^2 <\ell(\mathbf x) < 5id^2 + 4d^2$.
\end{itemize}
By construction of $C$ we have $4d^2 < \ell(\mathbf x) < 5d^3+d^ 2$ for all $\mathbf x \in C$, and so these cases cover all points in $C$. To show that each case actually occurs, we will construct a line segment $L$ along which $\ell$ takes values covering all cases above. 

Let $p$ be the point whose first coordinate is $5d^2-d +\frac{1}{2}$, and whose remaining coordinates are all equal to $1+\frac{1}{d-1}$. Observe that the sum of the coordinates of $p$ are exactly $5d^2 +\frac{1}{2}$. Moreover, we have \begin{align*}
\ell(p)&= 5d^2 - d + \frac{1}{2} + \left(1+\frac{1}{d-1}\right)\sum_{i=2}^d i\\
&=  5d^2 - d + \frac{1}{2} + \left(1+\frac{1}{d-1}\right)\left(d-1 + \sum_{i=1}^{d-1} i\right)\\
&\le  5d^2 - d + \frac{1}{2} + \left(1+\frac{1}{d-1}\right)\left(d-1 + (d-1)^2\right)\\
&=  5d^2 - d + \frac{1}{2} +d^2\\
& < 5d^2 + d^2.  
\end{align*}
Symmetrically, let $q$ be a point whose last coordinate is $5d^2-d +\frac{1}{2}$, and whose remaining coordinates are equal to $1+\frac{1}{d-1}$. As with $p$, we see that the sum of coordinates of $q$ is exactly $5d^2 + \frac{1}{2}$. Furthermore, we may compute\begin{align*}
\ell(q) & = d\left(5d^2-d +\frac{1}{2}\right) + \sum_{i=1}^{d-1} i\left(1+\frac{1}{d-1}\right)\\
& = 5d^3 -d^2 + \frac{d}{2} + \left(1+\frac{1}{d-1}\right)\sum_{i=1}^{d-1} i\\
& \ge 5d^3 - d^2 + \frac{d}{2} + d\\
& > 5d^3 - d^2. 
\end{align*}
 Since $\ell$ is linear, we conclude that $\ell$ takes all real values between $\ell(p) < 5d^2+d^2$ and $\ell(q) > 5d^3-d^2$ on the line segment $L$. Moreover, every point on $L$ has sum of coordinates equal to $5d^2+\frac{1}{2}$, so $L$ is contained in $C$. In particular, there are points in $L\cap C$ for which $\ell$ attains a value covering each case in the bulleted list above. Thus all codewords of types (i)-(iv) arise along $L$, and no others arise from the various $U_{\beta_i}$ and $U_{\gamma_i}$. 
 
 We have determined that the various $U_{\alpha_i}$ and $U_{\delta_i}$ give rise to exactly the codewords of types (vi)-(vii) in isolation, while the various $U_{\beta_i}$ and $U_{\gamma_i}$ give rise to exactly those of types (i)-(iv) in isolation. We must now argue that when considered together, all these codewords remain, and the only new codewords that appear are exactly those of type (v) in Definition \ref{def:(2,d,d)}. 
 
 Considering the sets together, we do not lose any codewords. Those of types (vi)-(vii) arise outside of $C$, which contains the various $U_{\beta_i}$ and $U_{\gamma_i}$. Those of types (i)-(iv) arise along $L$, and every point in $L$ has all coordinates larger than one, so $L$ does not meet any $U_{\alpha_i}$ or $U_{\delta_i}$. To see that only codewords of type (v) arise when considering the sets together, it will suffice to show that $U_{\alpha_i}\cap C$ is nonempty, and is contained in the region where the codeword $\{\beta_i,\gamma_i\}$ arises (recall that $U_{\alpha_i}$ and $U_{\delta_i}$ are identical outside of the unit hypercube, which does not meet $C$). That is, it will suffice to show that $5id^2-d^2 \le \ell(\mathbf x) \le 5id^2+d^2$ for all $\mathbf x \in U_{\alpha_i}\cap C$, and to find an example of one such $\mathbf x$. 
 
 The points in $\mathbf x$ in $U_{\alpha_i} \cap C$ are exactly those which satisfy the following $d+1$ conditions:\[
 5d^2 < \sum_{i\in[d]} x_i < 5d^2 + 1, \quad\quad
 0 < x_j < 1 \text{ for $j \neq i$, and} \quad\quad  x_i > 0.\]
Note that $\ell(\mathbf x)$ is smallest when the early coordinates of $\mathbf x$ are larger and the overall sum of the coordinates is smallest. Thus for $\mathbf x \in U_{\alpha_i}\cap C$, the value of $\ell(\mathbf x)$  is bounded below by \begin{align*}
 \sum_{j=1}^{i}j + i \big(5d^2-i \big) & = 5id^2 -i^2 + \sum_{j=1}^{i}j  \ge 5id^2 - d^2. 
 \end{align*}
 On the other hand, $\ell(\mathbf x)$ is largest when the later coordinates of $\mathbf x$ are larger, and the overall sum of the coordinates is largest. Thus the value of $\ell(\mathbf x)$ on the region $U_{\alpha_i} \cap C$ is bounded above by \begin{align*}
 i(5d^2 + 1 - i) + \sum_{j=i}^{d} j & = 5id^2 + i - i^2 + \sum_{j=i}^d j < 5id^2 + d^2. 
 \end{align*}
 This shows that $U_{\alpha_i}\cap C$ is contained in the region where the codeword $\{\beta_i,\gamma_i\}$ arises, so the only codeword that could arise in our overall realization involving both $\alpha_i$ and some $\beta_j$ is exactly $\{\alpha_i,\beta_i,\gamma_i,\delta_i\}$. To see that this codeword actually does arise, consider the point whose $i$-th coordinate is $5id^2+ \frac{1}{2} - \varepsilon$, and all of whose other coordinates are $\frac{\varepsilon}{d-1}$. For a sufficiently small $\varepsilon > 0$, this point will lie in $C\cap U_{\alpha_i}$, and thus give rise to the codeword $\{\alpha_i,\beta_i,\gamma_i,\delta_i\}$. 
 
 So far we have shown that our collection is an open convex realization of $\C_{(2,d,d)}$. Let us finally argue that our realization is non-degenerate. It suffices to observe that replacing the sets in our realization with their closures does not change the realized code. The arguments above can be applied verbatim, provided that we swap any strict inequalities for non-strict inequalities, and vice versa. 
\end{proof}

\begin{example}\label{ex:(2,3,3)}
Let us consider the code $\C_{(2,d,d)}$ in the case $d=3$. We have \begin{align*}
\C_{(2,3,3)} =  &\{\beta_1\gamma_1, \, \, \alpha_1\beta_1\gamma_1\delta_1, \, \,  \beta_{12}\gamma_1, \, \, \beta_{12}, \, \,  \beta_{12}\gamma_2, \\
& \, \,\,  \beta_2\gamma_2, \, \,  \alpha_2\beta_2\gamma_2\delta_2, \, \, \beta_{23}\gamma_2, \, \,  \beta_{23}, \, \,  \beta_{23}\gamma_3, \\ 
& \,\,\, \beta_3\gamma_3,\,\, \alpha_3\beta_3\gamma_3\delta_3, \\
&\,\,\, \alpha_1\delta_1,\,\, \alpha_2\delta_2,\,\, \alpha_3\delta_3,\,\, \delta_{123}, \,\, \alpha_{123}\delta_{123},\, \, \emptyset\}.
\end{align*}
Figure \ref{fig:(2,3,3)} illustrates the construction used in Proposition \ref{prop:(2,d,d)nondeg} to obtain a non-degenerate realization of $\C_{(2,3,3)}$ in $\R^3$. Note that Figure \ref{fig:(2,3,3)} is only a sketch of our construction---we do not precisely illustrate the inequalities that define the set $C$, and the various $U_{\alpha_i}$ and $U_{\delta_i}$ would be thinner relative to $C$ in an exact illustration.
\begin{figure}
\[\includegraphics{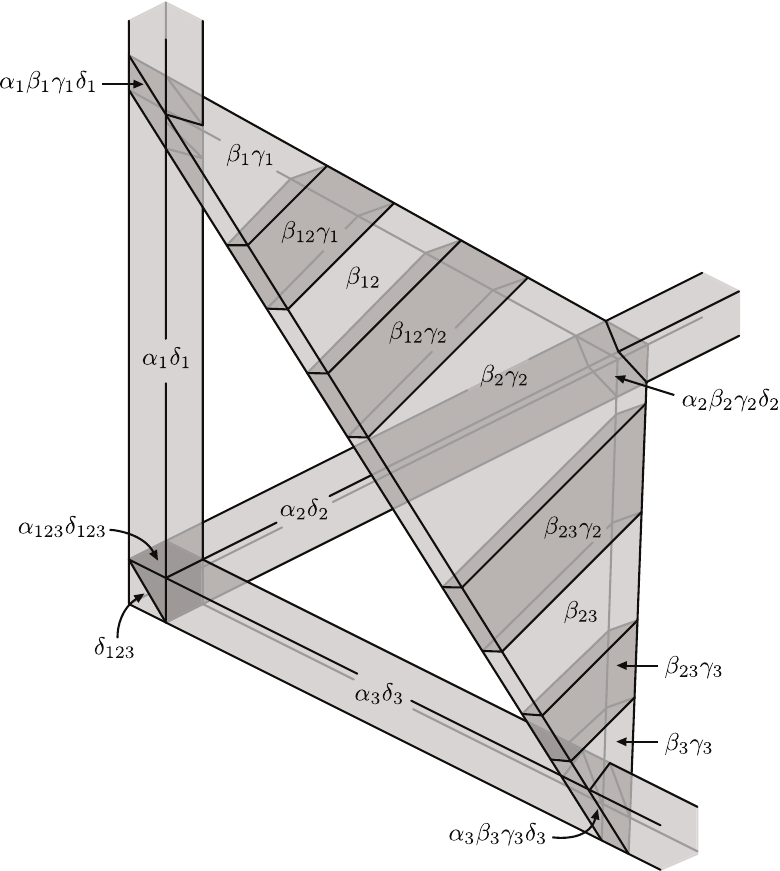}\]
\caption{A non-degenerate open realization of $\C_{(2,3,3)}$ in $\R^3$, with the regions that give rise to each codeword labeled.}\label{fig:(2,3,3)}
\end{figure}
\end{example}

So far, we have established appropriate upper bounds on the embedding dimensions of $\C_{(2,d,d)}$. We now move on to establish matching lower bounds. It will suffice to show that $\cdim(\C_{(2,d,d)}) = d$, which we do in Proposition \ref{prop:(2,d,d)closed}. Our proof requires two existing tools. The first tool is a notion of ``rigid structures" defined by \cite{rigidstructures}, which guarantees that a union of certain sets in a closed realization is convex---we do not state the definition of a rigid structure in full generality, but instead give a sufficient version of this result as a lemma below. The second tool we require is a code $\A_d$ from \cite{nonmonotonicity}. The relevant feature of this code is that if we add a certain codeword to it, the resulting code has closed embedding dimension equal to $d$---the codeword of type (viii) from Definition \ref{def:(2,d,d)} will be exactly the codeword that we need, up to a relabeling of the base set.

\begin{lemma}[Version of Lemma 4.21 from \cite{rigidstructures}]\label{lem:rigid}
Let $\X = \{X_1, X_2, \ldots, X_n\}$ be a closed convex realization of a code $\C$, and suppose that the nonempty codewords in $\C$ can be labeled $c_1, c_2, \ldots, c_k$ so that (i) $c_1\subset c_2 \supset c_3 \subset c_4 \supset \cdots \subset c_{k-1}\supset c_k$, (ii) no other containments occur between nonempty codewords, and (iii) $c_i\cap c_{i+1}\cap c_{i+2}$ is nonempty for all $i\in [k-2]$. Then the union $\bigcup_{i\in[n]} X_i$ is convex. 
\end{lemma}

\begin{proposition}\label{prop:(2,d,d)closed}
The code $\C_{(2,d,d)}$ has closed embedding dimension equal to $d$.
\end{proposition}
\begin{proof}
In Proposition \ref{prop:(2,d,d)nondeg} we showed that $\C_{(2,d,d)}$ has a non-degenerate realization in $\R^d$. Thus it will suffice to show that $\C_{(2,d,d)}$ does not have a closed realization in any dimension $d'<d$. Suppose for contradiction that we have a closed convex realization \[\X = \{X_{\alpha_1}, \ldots, X_{\alpha_d},X_{\beta_1}, \ldots, X_{\beta_d},X_{\gamma_1}, \ldots, X_{\gamma_d},X_{\delta_1}, \ldots, X_{\delta_d}\}\] of $\C_{(2,d,d)}$ in $\R^{d'}$ where $d'< d$. 

Consider the code that arises only from the various $X_{\beta_i}$ and $X_{\gamma_i}$. The nonempty codewords in this code will be exactly the codewords of types (i)-(iv) in Definition \ref{def:(2,d,d)}. Observe that we may order these codewords sequentially so that we have the containments \begin{align*}
\beta_1\gamma_1 \,\, \subset\, \, \beta_{12}\gamma_1 \,\, \supset \,\, \beta_{12}\,\,\subset\,\, \beta_{12}\gamma_2\,\,\supset \,\, \beta_2\gamma_2 \,\,\subset \,\, \beta_{23}\gamma_2\,\,\supset\,\, \cdots \\ \cdots \, \, \subset \beta_{d-1}\beta_d \gamma_{d-1} \supset \beta_{d-1}\beta_d\, \, \subset \beta_{d-1}\beta_d\gamma_d \, \, \supset \beta_d\gamma_d. 
\end{align*}
Moreover, no other containment relations exist between these codewords, and the intersection of any three consecutive codewords is nonempty (in particular, the intersection will contain some $\beta_i$). Thus by Lemma \ref{lem:rigid}, the union of all $X_{\beta_i}$ and $X_{\gamma_i}$ is a closed convex set. Let us call this union $X_{d+1}$. 

Now for $i\in[d]$, define $X_i = X_{\alpha_i}$ and $X_{\overline{i}} = X_{\delta_i}$. The code realized by the collection $\X' = \{X_1, X_2, \ldots, X_{d+1}, X_{\overline 1}, X_{\overline 2}, \ldots, X_{\overline d}\}$ will be exactly $\A_d \cup \{\{\overline 1, \overline 2, \ldots, \overline d\}\}$, where $\A_d$ is the code of \cite[Definition 3.6]{nonmonotonicity}. However, \cite[Theorem 3.7]{nonmonotonicity} states that the closed embedding dimension of $\A_d \cup \{\{\overline 1, \overline 2, \ldots, \overline d\}\}$ is exactly $d$. Since the realization $\X'$ lies in $\R^{d'}$ where $d'<d$, we have reached a contradiction. Thus $\C_{(2,d,d)}$ has closed embedding dimension equal to $d$. 
\end{proof}

We have established all the necessary constructions and results to exactly characterize the embedding dimensions of $\C_{(2,d,d)}$. We compile and summarize these results in the theorem below.

\begin{theorem}\label{thm:(2,d,d)}
The code $\C_{(2,d,d)}$ of Definition \ref{def:(2,d,d)} has embedding dimension vector equal to $(2,d,d)$. 
\end{theorem}
\begin{proof}
In Proposition \ref{prop:(2,d,d)open} we established that $\odim(\C_{(2,d,d)}) \le 2$, and in Proposition \ref{prop:(2,d,d)nondeg} we showed that $\nddim(\C_{(2,d,d)}) \le d$. Proposition \ref{prop:(2,d,d)closed} showed that $\cdim(\C_{(2,d,d)}) = d$, which implies that the non-degenerate embedding dimension is also equal to $d$. We cannot have $\odim(\C_{(2,d,d)}) < 2$ since Theorem \ref{thm:1dim} would imply that the closed embedding dimension is less than $d$. Thus the open embedding dimension is exactly $2$, and the result follows. 
\end{proof}

\section{Constructing the Codes $\C_{(2,2,d)}$ with $d< \infty$}\label{sec:(2,2,d)}

The code $\C_{(2,2,d)}$ is closely related to the code $\C_{(2,d,d)}$ which we defined and analyzed in the previous section. In fact, $\C_{(2,2,d)}$ is simply the result of deleting the base set elements $\delta_{[d]}$ from $\C_{(2,d,d)}$. It turns out this is enough to lower the closed embedding dimension from $d$ to $2$, without changing the other embedding dimensions. 

\begin{definition}\label{def:(2,2,d)}
Let $2\le d < \infty$, and define $\C_{(2,2,d)}$ to be the code on the base set $\alpha_{[d]}\cup \beta_{[d]}\cup \gamma_{[d]}$ which has following nonempty codewords:\begin{itemize}
\item[(i)] $\{\beta_i, \gamma_i\}$ for $i\in [d]$,
\item[(ii)] $\{\beta_i,\beta_{i+1}, \gamma_i\}$ for $i\in[d-1]$,
\item[(iii)] $\{\beta_i,\beta_{i+1}\}$ for $i\in[d-1]$,
\item[(iv)] $\{\beta_i,\beta_{i+1}, \gamma_{i+1}\}$ for $i\in[d-1]$,
\item[(v)] $\{\alpha_i, \beta_i,\gamma_i \}$ for $i\in [d]$,
\item[(vi)] $\{\alpha_i\}$ for $i\in [d]$,
\item[(vii)] $\alpha_{[d]}$.
\end{itemize}
\end{definition}

The close relationship between $\C_{(2,2,d)}$ and $\C_{(2,d,d)}$ greatly simplifies our analysis of $\C_{(2,2,d)}$. As a start, we have the following. 

\begin{proposition}\label{prop:(2,2,d)openandnondeg}
The code $\C_{(2,2,d)}$ has an open realization in $\R^2$ and a non-degenerate realization in $\R^d$.
\end{proposition}
\begin{proof}
Since $\C_{(2,2,d)}$ is the result of deleting the base set elements $\{\delta_1,\ldots, \delta_d\}$ from $\C_{(2,d,d)}$, any realization of $\C_{(2,d,d)}$ yields a realization of $\C_{(2,2,d)}$ by deleting the various $U_{\delta_i}$. In Proposition \ref{prop:(2,d,d)open} we constructed an open realization of $\C_{(2,d,d)}$ in $\R^2$, and in Proposition \ref{prop:(2,d,d)nondeg} we constructed a non-degenerate realization of $\C_{(2,d,d)}$ in $\R^d$. Since deleting sets in a realization preserves openness and non-degeneracy of the realization, these constructions give us an open realization of $\C_{(2,2,d)}$ in $\R^2$ and a non-degenerate realization of $\C_{(2,2,d)}$ in $\R^d$ as desired.
\end{proof}

It remains to argue that $\C_{(2,2,d)}$ has a closed realization in $\R^2$, but no non-degenerate realization in a dimension less than $d$. We start by constructing a closed realization. Example \ref{ex:(2,2,d)closed} illustrates this construction in the case $d=4$. 

\begin{proposition}\label{prop:(2,2,d)closed}
The code $\C_{(2,2,d)}$ has a closed realization in $\R^2$.
\end{proposition}
\begin{proof}
Informally, we may form a closed realization by arranging the various $X_{\beta_i}$ and $X_{\gamma_i}$ sequentially along the $x$-axis, and then letting the various $X_{\alpha_i}$ be triangles which meet at a common point above the $x$-axis and intersect the $x$-axis sequentially. Formally, let $C$ be the strip $\{(x,y)\in \R^2\mid -1\le y \le 0 \text{ and } 0\le x\le 4d-3\}$. Then for $i\in [d]$ we define \begin{align*}
X_{\beta_i} & = \{(x,y)\in C \mid  4i-7 \le x \le 4i\}, \text{and}\\
X_{\gamma_i}  & =\{(x,y)\in C \mid 4i-5 \le x \le 4i-2\}.
\end{align*}  
Let $p = (0,d)$, and for $i\in[d]$ let $q_i = (4i-3.75, 0)$ and $r_i = (4i-3.25, 0)$. Then define \[
X_{\alpha_i} = \conv\{p, q_i, r_i\}. 
\]
We claim that this yields a closed realization of $\C_{(2,2,d)}$. Observe that the $X_{\alpha_i}$ are triangles which only meet at $p$, so the codewords they give rise to in isolation are simply $\{\alpha_i\}$ for $i\in[d]$ and $\alpha_{[d]}$, the latter arising only at $p$. These are exactly the codewords of types (vi) and (vii) in Definition \ref{def:(2,2,d)}. 

The codewords that arise from the various $U_{\beta_i}$ and $U_{\gamma_i}$ in isolation are completely determined by the $x$ coordinates of points in $C$. Indeed, if $(x,y)$ is a point in $C$ then the codeword arising at this point is\begin{itemize}
\item $\{\beta_i,\gamma_i\}$ if and only if $4i-4 < x < 4i-3$,
\item $\{\beta_i, \beta_{i+1}, \gamma_i\}$ if and only if $4i-3 \le x \le 4i-2$,
\item $\{\beta_i, \beta_{i+1}\}$ if and only if $4i-2 < x < 4i-1$, and
\item $\{\beta_i,\beta_{i+1}, \gamma_{i+1}\}$ if and only if $4i-1\le x \le 4i$. 
\end{itemize}
These cases partition all points in $C$, and all such cases occur by construction of $C$. Moreover, these are exactly the codewords of types (i)-(iv) in Definition \ref{def:(2,2,d)}. Finally, note that by choice of the points $q_i$ and $r_i$, the triangle $U_{\alpha_i}$ only meets $C$ where the codeword $\{\beta_i,\gamma_i\}$ arises. This yields exactly the codewords of type (v), and so we have indeed constructed a closed realization of $\C_{(2,2,d)}$ as desired.
\end{proof}

\begin{example}\label{ex:(2,2,d)closed}
Figure \ref{fig:(2,2,d)closed} shows the construction used in Proposition \ref{prop:(2,2,d)closed} to obtain a closed realization of the code $\C_{(2,2,d)}$ in $\R^2$ in the case $d=4$. In this case, we have \begin{align*}
\C_{(2,2,4)} = & \{\beta_1\gamma_1, \, \, \alpha_1\beta_1\gamma_1, \, \,  \beta_{12}\gamma_1, \, \, \beta_{12}, \, \,  \beta_{12}\gamma_2, \\
& \, \,\, \beta_2\gamma_2, \, \,  \alpha_2\beta_2\gamma_2, \, \, \beta_{23}\gamma_2, \, \,  \beta_{23}, \, \,  \beta_{23}\gamma_3, \\
& \, \, \,
\beta_3\gamma_3, \, \,  \alpha_3\beta_3\gamma_3, \, \,  \beta_{34}\gamma_3, \, \,  \beta_{34}, \, \,  \beta_{34}\gamma_4,\\
&\,\,\, \beta_4\gamma_4, \,\, \alpha_4\beta_4\gamma_4, \\
& \, \,\, \alpha_1, \, \,  \alpha_2, \, \,  \alpha_3, \, \,  \alpha_4,  \, \,  \alpha_{1234},\, \, \emptyset\}.
\end{align*}
\begin{figure}[h]
\[\includegraphics{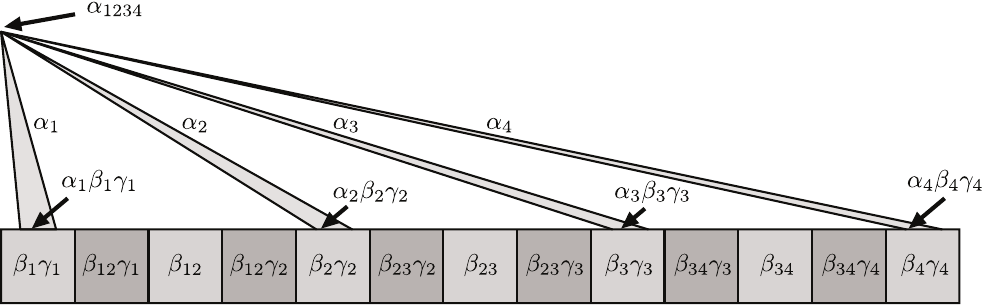}\]
\caption{A closed realization of $\C_{(2,2,4)}$ in $\R^2$.}\label{fig:(2,2,d)closed}
\end{figure}
\end{example}

We are now ready to prove that $\nddim(\C_{(2,2,d)}) = d$. Our proof proceeds similarly to the proof of Proposition \ref{prop:(2,d,d)closed}---namely, it relies on the rigid structure result in Lemma \ref{lem:rigid}, and on the characterization of the embedding dimensions of an existing family of codes. 

\begin{proposition}\label{prop:(2,2,d)nondeg}
The code $\C_{(2,2,d)}$ has non-degenerate embedding dimension equal to $d$.
\end{proposition}
\begin{proof}
In Proposition \ref{prop:(2,2,d)openandnondeg} we argued that $\nddim(\C_{(2,2,d)} \le d$, so it will suffice to argue that there is no non-degenerate (open or closed) realization of $\C_{(2,2,d)}$ in $\R^{d'}$ with $d'<d$. Suppose for contradiction that we have a closed non-degenerate realization \[
\X = \{X_{\alpha_1}, \ldots, X_{\alpha_d}, X_{\beta_1}, \ldots, X_{\beta_d},X_{\gamma_1},\ldots, X_{\gamma_d}\}
\]
of $\C_{(2,2,d)}$ in a dimension $d'<d$. As in the proof of Proposition \ref{prop:(2,d,d)closed}, the codewords of types (i)-(iv) satisfy the conditions of Lemma \ref{lem:rigid}, and so the union of all $X_{\beta_i}$ and $X_{\gamma_i}$ is a closed convex set. Let us call this set $X_{d+1}$, and for $i\in[d]$ define $X_i = X_{\alpha_i}$. 

Now, we claim that non-degeneracy of the realization $\X$ guarantees non-degeneracy of the realization $\X' = \{X_1, X_2, \ldots, X_{d+1}\}$. First observe that the sets $\{X_1, \ldots, X_d\}$ are non-degenerate in isolation, and realize the code $\{[d], \{1\}, \{2\}, \ldots, \{d\}, \emptyset \}$. The set $X_{d+1}$ is full-dimensional, and since the codewords of type (v) in Definition \ref{def:(2,2,d)} are the only ones that simultaneously contain some $\alpha_i$ and $\beta_j$ or $\gamma_j$, we conclude that $X_{d+1}$ intersects each other $X_i$ at a point common to both their interiors, while avoiding the common intersection of all $X_i$ for $i\in[d]$.

However, the analysis above tells us that $\code(\X') = \S_d$ (recall Definition \ref{def:Sd}). The code $\S_d$ has non-degenerate embedding dimension exactly $d$, and since $\X'$ is a realization in $\R^{d'}$ with $d'< d$, we have reached a contradiction. This proves the result. 
\end{proof}

We have now characterized the embedding dimension vector of $\C_{(2,2,d)}$, as summarized in the theorem below.

\begin{theorem}\label{thm:(2,2,d)}
The code $\C = \C_{(2,2,d)}$ of Definition \ref{def:(2,2,d)} has embedding dimension vector equal to $(2,2,d)$. 
\end{theorem}
\begin{proof}
In Proposition \ref{prop:(2,2,d)openandnondeg} we showed that the open and non-degenerate embedding dimensions of $\C_{(2,2,d)}$ were no larger than two and $d$, respectively. Proposition \ref{prop:(2,2,d)closed} showed that the closed embedding dimension was no more than two. In Proposition \ref{prop:(2,2,d)nondeg} we argued that the non-degenerate embedding dimension was exactly $d$, which then implies that the closed and open embedding dimensions are both exactly two---if they were smaller, so would be the non-degenerate embedding dimension.
\end{proof}

\section{Constructing the Codes $\C_{(\infty, 2, \infty)}$, $\C_{(2,\infty,\infty)}$, and $\C_{(2,2,\infty)}$}\label{sec:infty}

We now treat the three remaining cases, in which some embedding dimensions may be infinite. As we did in Section \ref{sec:(d,2,d)}, we draw on some existing examples in the literature which suffice---in fact, the codes $\C_{(\infty,2,\infty)}$ and $\C_{(2,\infty,\infty)}$ have already been defined an analyzed in \cite{obstructions} and \cite{openclosed} respectively. Our main contribution is the construction of the code $\C_{(2,2,\infty)}$ (see Theorem \ref{thm:(2,2,infty)}), which adds a rigid structure to the minimally non-convex code from \cite[Theorem 5.10]{morphisms}.

\begin{figure}[h]
\[\includegraphics{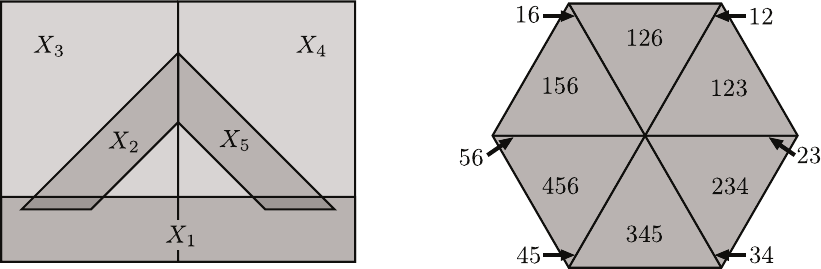}\]
\caption{Left: a closed realization of $\C_{(\infty, 2, \infty)}$. Right: an open realization of $\C_{(2,\infty,\infty)}$. In the open realization we have labeled the regions that give rise to each codeword---the various $U_i$ are open halves of the hexagon, rotated consecutively by $60$ degrees.}\label{fig:inftyexamples}
\end{figure}

\begin{proposition}\label{prop:(infty,2,infty)}
The code \[
\C_{(\infty, 2,\infty)} \od \{\mxl{2345}, \mxl{123}, \mxl{134}, \mxl{145}, 13, 14, 23, 34, 45, 3, 4, \emptyset\}
\]
which appears in \cite[Theorem 3.1]{obstructions} has embedding dimension vector $(\infty, 2, \infty)$. 
\end{proposition}
\begin{proof}
\cite[Theorem 3.1]{obstructions} states that this code does not have an open convex realization in any dimension, and so its open and non-degenerate embedding dimensions are both $\infty$. On the other hand, \cite[Figure 2.2(a)]{openclosed} provides a closed realization of this code in $\R^2$. 
\end{proof}

\begin{proposition}\label{prop:(infty,2,infty)}
The code \[
\C_{(2,\infty,\infty)} \od \{\mxl{123}, \mxl{126}, \mxl{156}, \mxl{234}, \mxl{345}, \mxl{456}, 12, 16, 23, 34, 45, 56, \emptyset\}
\]
which appears in \cite[Section 2.3]{openclosed} has embedding dimension vector $(2, \infty, \infty)$. 
\end{proposition}
\begin{proof}
\cite[Lemma 2.9]{openclosed} states that this code does not have a closed convex realization in any dimension, and so its closed and non-degenerate embedding dimensions are both $\infty$. However, \cite[Figure 2.1(a)]{openclosed} provides an open realization of this code in $\R^2$. 
\end{proof}

Figure \ref{fig:inftyexamples} duplicates \cite[Figure 1.7]{amziphdthesis}, illustrating a closed realization of $\C_{(\infty, 2, \infty)}$ and an open realization of $\C_{(2,\infty,\infty)}$ in $\R^2$. Note that both realizations are degenerate. For example, on the left $X_2$ and $X_3$ only intersect in a 1-dimensional segment, and on the right $U_1$ and $U_4$ are disjoint but share boundary points. In Theorem \ref{thm:(2,2,infty)}, we conclude by constructing and analyzing the code $\C_{(2,2,\infty)}$. 

\begin{theorem}\label{thm:(2,2,infty)}
The code \[
\C_{(2,2,\infty)} \od \{\mxl{123}, \mxl{145}, \mxl{2456}, \mxl{2467}, \mxl{389}, \mxl{678}, \mxl{689}, 246, 45, 67, 68, 89, 1, 2, 3, \emptyset\}
\]
has embedding dimension vector $(2,2,\infty)$.
\end{theorem}
\begin{proof}

\begin{figure}[h]
\[\includegraphics{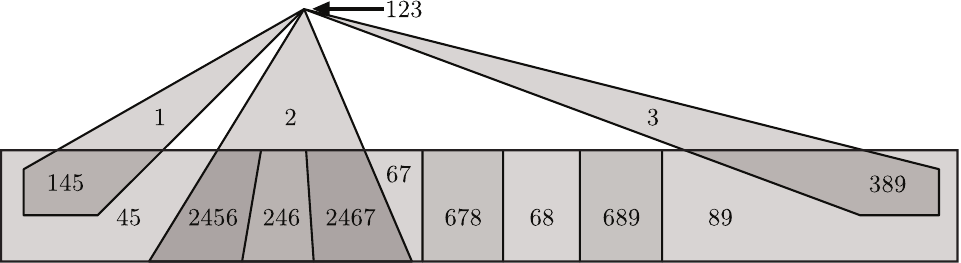}\]\vspace{-2em}
\caption{A closed realization of $\C_{(2,2,\infty)}$ in $\R^2$.}\label{fig:(2,2,infty)closed}
\end{figure}

\begin{figure}[h]
\[\includegraphics{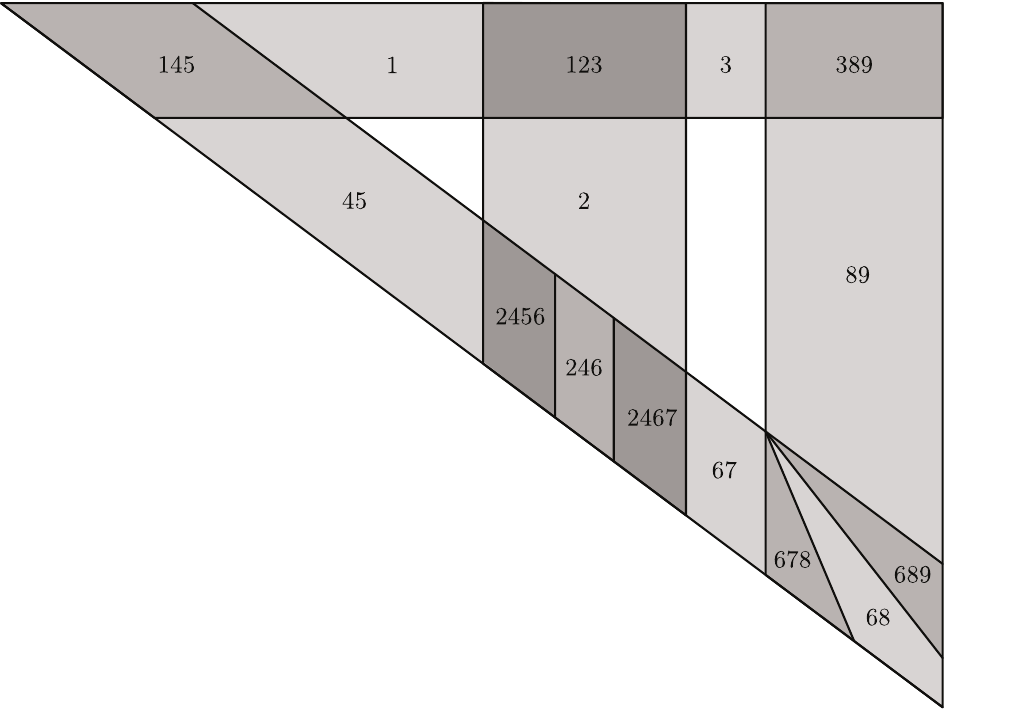}\]\vspace{-4em}
\caption{An open realization of $\C_{(2,2,\infty)}$ in $\R^2$.}\label{fig:(2,2,infty)open}
\end{figure}

Figures \ref{fig:(2,2,infty)closed} and \ref{fig:(2,2,infty)open} show closed and open realizations of $\C_{(2,2,\infty)}$ in $\R^2$. It remains to show that no non-degenerate realization of this code exists. Suppose for contradiction that we have a closed non-degenerate realization $\X = \{X_1, X_2, \ldots, X_9\}$ in $\R^d$. We may assume without loss of generality that the various $X_i$ are compact.

 Let us first examine the sets $\{X_4, X_5, \ldots, X_9\}$ in isolation. The codewords that arise from this collection will be $\{45, 456, 46, 467, 67, 678, 68, 689, 89,\emptyset\}$. These codewords satisfy the conditions of Lemma \ref{lem:rigid}---namely, we have the containments \[
45\, \subset \,  456\, \supset\, 46\,\subset\, 467\,\supset\, 67\,\subset\, 678\,\supset\, 68\,\subset\, 689\,\supset\, 89,
\]
no other containments occur, and the intersection of any three consecutive codewords is nonempty. Thus the union $X_4\cup X_5\cup \cdots \cup X_9$ is convex. Since our realization is non-degenerate, the interior of this union is a nonempty convex open set---let us call this interior $U$.

Now, let $p$ be a point in the interior of $X_{123}$, let $q$ be a point in the interior of $X_{145}$, and let $r$ be a point in the interior of $X_{389}$. Observe that $q$ and $r$ both lie in $U$, so the line segment $\overline{qr}$ is contained in $U$. Thus the consecutive codewords that appear along $\overline{qr}$ must all contain some index between $4$ and $9$. In fact, consecutive codewords that appear along this line segment must contain one another (see \cite[Lemma 2.1]{orderforcing}). The only possible sequence of codewords that can arise along $L$ is therefore \[
145 \quad 45 \quad 2456 \quad 246 \quad 2467\quad 67 \quad 678 \quad 68 \quad 689 \quad 89 \quad 389. 
\]
In particular, $L$ passes through $X_2$ in addition to the interiors of $X_1$ and $X_3$. By possibly perturbing $q$ and $r$ by a small distance, we may assume that $\overline{qr}$ passes through the interior of $X_2$. 

Let $A$ be the affine span of $p$, $q$, and $r$. We may assume that these points are in general position so that $A$ has dimension exactly two. Define $U_i = \interior(X_i)\cap A$ for $i\in [3]$, and $U_4 = U\cap A$. The sets $U_1$, $U_2$, $U_3$, and $U_4$ may then be regarded as convex open sets in $A\cong \R^2$.

We claim that the code realized by $\U = \{U_1, U_2, U_3, U_4\}$ is exactly $\S_3$ (recall Definition \ref{def:Sd}). Our choice of $p$ guarantees that $U_1\cap U_2\cap U_3$ is nonempty, while the line segment $\overline{qr}\subseteq U_4$ guarantees that $U_4$ intersects each of $U_1, U_2,$ and $U_3$. The set $U_4$ does not meet $U_{123}$ since $U_4$ is contained in $U$ which does not meet $X_{123}$, which in turn contains $U_{123}$. Finally, codewords containing any of $12, 13,$ or $23$ do not appear in this realization---if they did, then non-degeneracy of $\X$ would imply that there was some codeword appearing in the original realization which contained one of these. No such codeword exists in the original realization, so $\U$ is an open convex realization of $\S_3$ in $\R^2$. This contradicts the fact that $\odim(\S_3) = 3$ (recall Proposition \ref{prop:(d,2,d)}). Thus $\C_{(2,2,\infty)}$ does not have a non-degenerate realization in any dimension.\end{proof}

\section{Conclusion}
We have constructed several families of codes and characterized their embedding dimension vectors. In combination, these families guarantee that every vector $(a,b,c)$ with $2\le a,b,c\le \infty$ and $\max\{a,b\} \le c$ arises as the embedding dimension vector of some code (Theorem \ref{thm:main}). Moreover, such vectors are \emph{exactly} those that can arise as embedding dimension vectors, with the exception of the vectors $(0,0,0)$ and $(1,1,1)$. Although our results required careful and sometimes lengthy proofs, our arguments were primarily based on existing tools in the convex neural code literature: sunflowers of convex open sets (recall Section \ref{sec:(d,2,d)}, which restates results of \cite{embeddingphenomena}), and rigid structures of closed convex sets (recall Lemma \ref{lem:rigid}, which restates results of \cite{rigidstructures}). Our contribution was to find combinations of these tools which yielded the correct embedding dimension vectors, and carry out the necessary analysis to characterize these vectors.

Although we have settled the question of which vectors can arise as embedding dimension vectors, it is still very difficult to bound the embedding dimensions of an arbitrary code $\C\subseteq 2^{[n]}$. In fact, it is even an open question whether or not there exists an algorithm which can decide the open, closed, or non-degenerate embedding dimension of a code.

A further area of study which we did not explore in this work was the relationship between the size of the base set of a code and its embedding dimensions. Such a line of investigation would help characterize how ``efficiently" codes can capture the dimension of a space in which they are realized. We thus ask the following:

\begin{question}\label{q:max}
Among all codes $\C$ with base set $[n]$, what is the maximum finite open (respectively, closed or non-degenerate) embedding dimension that arises? Which codes achieve this maximum?
\end{question}

As a start, we conjecture that each additional base set element yields a strict increase in the maximum embedding dimension. 

\begin{conjecture}\label{conj:strict}
The maximum described in Question \ref{q:max} is a strictly increasing function of $n$. That is, if $\C\subseteq 2^{[n]}$ has maximum open embedding dimension among all codes on $[n]$, while $\D\subseteq 2^{[n+1]}$ has maximum open embedding dimension among all codes on $[n+1]$, then $\odim(\D) > \odim(\C)$. Moreover, the analogous result should hold for closed and non-degenerate embedding dimensions. 
\end{conjecture}

Rather than stratifying codes by the size of their base sets, and then asking for the maximum embedding dimension in each strata, one could take the reverse perspective: stratify by embedding dimension, and then ask for the smallest base set size. We formalize this point of view below. Note that this is not simply a reformulation of Question \ref{q:max}, though these two questions do provide bounds for one another. 

\begin{question}\label{q:min}
Among all codes $\C$ with open (respectively, closed or non-degenerate) embedding dimension equal to $d$, what is the minimum base set size that arises? Which codes achieve this minimum?
\end{question}

Results of \cite{embeddingphenomena} imply that the largest finite open embedding dimension among codes $\C\subseteq 2^{[n]}$ can be as large as $\binom{n-1}{\lfloor (n-1)/2\rfloor}$---in particular, it can be larger than $n$. This implies that the minimum in Question \ref{q:min} is \emph{not} a strictly increasing function of $d$, so we cannot make an analogous conjecture to Conjecture \ref{conj:strict} in the case of open embedding dimension. However, there are not yet known examples where $\cdim(\C) > n$ and $\C\subseteq 2^{[n]}$. Nevertheless, we conjecture that such codes exist, and finding such examples would be a good starting point for work on Question \ref{q:min}.

\begin{conjecture}
There exists a code $\C\subseteq 2^{[n]}$ such that $n < \cdim(\C) < \infty$. 
\end{conjecture}

\section*{Acknowledgements}

We thank Patrick Chan, Katherine Johnston, Joseph Lent, Alexander Ruys de Perez, and Anne Shiu for sharing early drafts of their work on rigid structures. We especially thank Anne Shiu for discussion on this topic, which helped formalize and streamline the presentation of our families of codes. We thank Florian Frick for asking questions which motivated us to formulate Questions \ref{q:max} and \ref{q:min}. We are also grateful to the anonymous referees for helpful feedback and suggestions. 

\bibliographystyle{plain}
\bibliography{references}

\end{document}